# A note on the polar decomposition in metric spaces

Zhirayr Avetisyan[*], Michael Ruzhansky[†]

September 2, 2023

## Abstract

The analogue of polar coordinates in the Euclidean space, a polar decomposition in a metric space, if well-defined, can be very useful in dealing with integrals with respect to a sufficiently regular measure. In this note we handle the technical details associated with such polar decompositions.

## 1   Introduction

From a bird's perspective, in dealing with integrals of measurable functions with respect to infinite measures, it is often technically very convenient to view the full integral as a double integral, where the "infiniteness" of the measure falls completely onto the outer integral, e.g.,

$$\int_{\mathcal{M}} f(x)d\mu(x) = \int_{[0,R)} \int_{\mathcal{S}_r} f(x)d\omega_r(x)d\nu(r), \quad R \in (0,+\infty].$$

This assumes a measure space bundle $(\mathcal{M},\mu) \xrightarrow{\rho} ([0,R),\nu)$, $\mathcal{S}_r = \rho^{-1}(\{r\})$, and a disintegration

$$\mu(\cdot) = \int_{[0,R)} \omega_r(\cdot) d\nu(r),$$

with $\omega_r$ bounded and concentrated on $\mathcal{S}_r$ for $\nu$-a.e. $r \in [0,R)$. If $\mathcal{M}$ can be seen as a product $[0,R) \times \mathcal{S}$, i.e., there is a Borel almost isomorphism $([0,R) \times \mathcal{S}, \tilde{\mu}) \xrightarrow{\Psi} \mathcal{M}$ such that $\rho \circ \Psi = \pi_1$ (here $\pi_1(r,\varphi) = r$), then new variables $\mathcal{M} \ni x \mapsto (r,\varphi) \in [0,R) \times \mathcal{S}$ can be introduced, so that

$$\int_{\mathcal{M}} f(x)d\mu(x) = \int_{[0,R)} \int_{\mathcal{S}} \tilde{f}(r,\varphi)d\tilde{\omega}_r(\varphi)d\nu(r), \quad \tilde{f}(r,\varphi) = f \circ \Psi(r,\varphi), \quad d\tilde{\omega}_r(\varphi) = d\omega_r \circ \Psi(r,\varphi).$$

The choice of the function $\rho$ is pretty arbitrary at this level, but in metric spaces there is a natural candidate - the distance from a fixed point $x_0 \in \mathcal{M}$ of origin. If $\mathcal{M} = \mathbb{R}^n$ and $d\mu(x) = \dot{\mu}(x)dx$, $\dot{\mu} \in C^\infty(\mathbb{R}^n)$ a smooth positive density, then this simply amounts to switching to polar coordinates,

$$\int_{\mathbb{R}^n} f(x)d\mu(x) = \int_0^\infty \int_{\mathcal{S}_r} f(r,\varphi)\omega(r,\varphi)d\varphi dr, \quad \forall f \in L^1_{\mathrm{loc}}(\mathbb{R}^n),$$

$$\mathcal{S}_r \doteq \{x \in \mathbb{R}^n \mid |x| = r\}, \quad \omega(r,\varphi) = \frac{d\mu(x(r,\varphi))}{d\varphi dr}.$$

---

[*]Department of Mathematics: Analysis, Logic and Discrete Mathematics, Ghent University, Belgium, zhirayr.avetisyan@ugent.be

[†]Department of Mathematics: Analysis, Logic and Discrete Mathematics, Ghent University, Belgium, and School of Mathematical Sciences, Queen Mary University of London, United Kingdom. m.ruzhansky@ugent.be

More generally, let $\mathcal{M} = (\mathcal{M}, d)$ be a metric space and $\mu$ a measure on the Borel $\sigma$-algebra $\Sigma$ of $\mathcal{M}$. All measures in this paper will be assumed positive, $\mu(\cdot) \geq 0$, and $\sigma$-additive, i.e.,

$$(\forall \{A_k\}_{k=1}^\infty \subset \Sigma) \quad (\forall k, l \in \mathbb{N})\, k \neq l \Rightarrow A_k \cap A_l = \emptyset \quad \Rightarrow \quad \mu\left(\bigcup_{k=1}^\infty A_k\right) = \sum_{k=1}^\infty \mu(A_k).$$

We will call every such measure a Borel measure, although some sources require that Borel measures have additional neat properties. Fix a point $x_0 \in \mathcal{M}$ and denote by $\rho_{x_0} \in C(\mathcal{M}, [0, R))$ the function given by

$$\rho_{x_0}(x) = d(x, x_0), \quad \forall x \in \mathcal{M}. \tag{1}$$

We will speak of spheres $\mathcal{S}_r(x_0)$ and balls $\mathcal{B}_r(x_0)$ of radius $r \in [0, R)$ centred at $x_0$,

$$\mathcal{S}_r(x_0) \doteq \rho_{x_0}^{-1}(\{r\}), \quad \mathcal{B}_r(x_0) \doteq \rho_{x_0}^{-1}([0, r)) = \bigcup_{0 \leq s < r} \mathcal{S}_s(x_0), \quad \forall r \in [0, R).$$

The Borel measure $\nu_{x_0}$ on $[0, R)$ is defined by

$$\nu_{x_0} \doteq \mu \circ \rho_{x_0}^{-1}, \tag{2}$$

so that

$$\nu_{x_0}([0, r)) = \mu\left(\mathcal{B}_r(x_0)\right), \quad \forall r \in [0, R).$$

**Definition 1** *We will say that the metric measure space $(\mathcal{M}, d, \mu)$ admits a polar decomposition at the point $x_0 \in \mathcal{M}$ if there exists a field $r \mapsto \omega_r$ of bounded (i.e., finite) Borel measures $\omega_r$ on $\mathcal{M}$ for $\nu_{x_0}$-almost all $r \in [0, R)$ such that for every measurable function $f : \mathcal{M} \to \mathbb{C}$, for which the integral*

$$\int_\mathcal{M} f(x)d\mu(x)$$

*makes sense (finite or infinite), the following conditions hold:*

1. *The map*

$$[0, R) \ni r \mapsto \int_\mathcal{M} f(x)d\omega_r(x)$$

   *is Borel measurable.*

2. *The measure $\omega_r$ is concentrated on $\mathcal{S}_r(x_0)$ for $\nu_{x_0}$-a.e $r \in [0, R)$, i.e.,*

$$\int_\mathcal{M} f(x)d\omega_r(x) = \int_{\mathcal{S}_r(x_0)} f(x)d\omega_r(x).$$

3. We have
$$\int_{\mathcal{M}} f(x)d\mu(x) = \int_{[0,R)} \int_{\mathcal{S}_r(x_0)} f(x)d\omega_r(x)d\nu_{x_0}(r)$$

**Remark 1** *If the measure $\nu_{x_0}$ is absolutely continuous with respect to the Lebesgue measure on $[0, R)$ then there exists a Borel measurable, non-negative, locally integrable function $\dot{\nu}_{x_0}$ such that $d\nu_{x_0}(r) = \dot{\nu}_{x_0}(r)dr$. In that case, the factor $\dot{\nu}_{x_0}$ can be incorporated into the measures $\omega_r$, so that part 3. of Definition 1 becomes*

$$\int_{\mathcal{M}} f(x)d\mu(x) = \int_0^R \int_{\mathcal{S}_r(x_0)} f(x)d\omega_r(x)dr.$$

*This corresponds to the definition (1.2) of polar decomposition in [RuVe19], so that in this case the further integration-by-parts techniques as used in that work are applicable.*

From a measure theory perspective, a polar decomposition is merely an instance of measure disintegration, and our job in this paper is to study sufficient conditions for a measure disintegration of this kind in a metric space. Whether or not $\nu_{x_0}$ is absolutely continuous with respect to the Lebesgue measure is unrelated to the disintegration and is a reflection of the geometric shape of $\mathcal{M}$ and the distribution of $\mu$.

## 2 A discussion of assumptions

**Finite measure of balls.** A polar decomposition as in Definition 1 is useful when the condition 3. is satisfied non-trivially. Let us write it in the measure form,

$$\mu(\cdot) = \int_{[0,R)} \int_{\mathcal{S}_r(x_0)} \omega_r(\cdot)d\nu_{x_0}(r).$$

If the measure $\mu$ is bounded then everything works well, but the point of introducing a polar decomposition is often to deal with an infinite measure. If the measure $\mu$ is not locally finite then certainly $C(\mathcal{M}) \not\subset L^1_{\text{loc}}(\mathcal{M}, \mu)$, and one is led to consider very thin spaces of integrable functions. Arguably, there is little use of such measures in functional analysis, and assuming the measure $\mu$ locally finite is not a serious loss of generality.

If $\mu$ is locally finite at $x_0$ then let

$$R_* \doteq \inf \{r \in [0, R) | \quad \mu(\mathcal{B}_r(x_0)) = +\infty\} \in (0, +\infty], \tag{3}$$

with the convention that $\inf \emptyset = R$. If $R_* < R$ then for every measurable set $A$ with $\rho_{x_0}(A) \subset [R_*, R)$ we have $\mu(A) = +\infty$, which makes the formula 3. in Definition 1 useless for $r \geq R_*$. As far as a polar decomposition around the fixed point $x_0$ is concerned, one can consider the subspace $\mathcal{M}' \doteq \mathcal{B}_{R^*}(x_0)$ instead of $\mathcal{M}$ and set $R \doteq R_*$. This way all balls $\mathcal{B}_r(x_0)$ for $r \in [0, R)$ will have finite volume. In particular, $\mu$ will be locally finite on $\mathcal{M}$.

**Separability and completeness of** $(\mathcal{M}, d)$**.** The separability of the metric space $(\mathcal{M}, d)$, and thus of its Borel $\sigma$-algebra and all Lebesgue functions on it, is of utmost importance for most functional analytical considerations. It is therefore not a serious restriction to assume $(\mathcal{M}, d)$ separable. However, we will see in Theorem 2 below that polar decompositions can be obtained even without separability, at the cost of completeness of $\mu$ and further restrictions.

If the metric space $(\mathcal{M}, d)$ is not complete then closed balls may fail to be compact and may have infinite measure even if $\mu$ is locally finite and $\mathcal{M}$ locally compact. Indeed, take $\mathcal{M} = (0, 1)$ with the Euclidean metric and the measure $d\mu(x) = dx/x$. Assuming the completeness of $\mathcal{M}$ wards against pathological behaviour of locally finite Borel measures. However, our main results, Theorem 1 and Theorem 2, will not assume the completeness of $(\mathcal{M}, d)$.

On the other hand, if we agree to embrace separability and completeness then we are in the realm of Polish metric spaces $(\mathcal{M}, d)$. Borel measures on Polish spaces are very tame and provide convenient technical ease.

**$\sigma$-compactness of $\mathcal{M}$.** Under this rubric we will present two technical statements which can be seen as partial converses of each other. Together they will show that $\sigma$-compactness of $\mathcal{M}$ is not too far from being the right context for polar decompositions, although it is not strictly necessary. Recall that a topological space is called $\sigma$-compact if it is the union of countably many compact subspaces.

First we will show that if we are willing to assume that the metric space $(\mathcal{M}, d)$ is Polish, then we could instead take $\mathcal{M}$ as $\sigma$-compact.

**Proposition 1** *Every locally finite Borel measure on a Polish space is concentrated on a $\sigma$-compact subspace, which is the union of a countable locally finite family of countable disjoint unions of compact sets.*

**Proof:** Let $\mathcal{M}$ be a Polish space and $\mu$ a locally finite Borel measure. Let $\{U_x\}_{x \in \mathcal{M}}$ be a family of open subsets such that $\mu(U_x) < \infty$ for $\forall x \in \mathcal{M}$, by local finiteness of $\mu$. As a metrizable space, $\mathcal{M}$ is paracompact [Sto48], and thus there is a locally finite open refinement $\{V_\alpha\}_{\alpha \in A}$ of $\{U_x\}_{x \in \mathcal{M}}$. But $\mathcal{M}$ is also second countable and hence Lindelöf, so that $\{V_\alpha\}_{\alpha \in A}$ has a countable locally finite subcover $\{V_k\}_{k=1}^{\infty}$.

Fix $k \in \mathbb{N}$. Clearly, $\exists x \in \mathcal{M}$ s.t. $V_k \subset U_x$, and therefore $\mu(V_k) \leq \mu(U_x) < \infty$. As an open subset of the Souslin space $\mathcal{M}$, $V_k$ is itself a Souslin space (Lemma 6.6.5 in [Bog07.II]), and the restriction of $\mu$ to $V_k$ is a bounded Borel measure. By Theorem 7.4.3 in [Bog07.II], $\mu$ is tight on $V_k$, i.e.,

$$(\forall \epsilon > 0)(\exists K_\epsilon \Subset V_k \text{ compact}) \, \mu(V_k \setminus K_\epsilon) < \epsilon.$$

Let $X_1 \doteq V_k$ and $K_1 \Subset X_1$ compact be chosen as above with $\mu(X_1 \setminus K_1) < 1$. Then $X_2 \doteq X_1 \setminus K_1$ is an open subset of a Souslin space, and thus a Souslin subspace, and $\mu$ restricted to $X_2$ gives a bounded Borel measure, which is again tight by the same theorem. Continuing this process inductively, on the step $n \in \mathbb{N}$ we choose $K_n \Subset X_n$ compact such that $\mu(X_n \setminus K_n) < \frac{1}{n}$. Denote

$$Y_k \doteq \bigcup_{n=1}^{\infty} K_n \subset X.$$

Then
$$\mu(V_k \setminus Y_k) \leq \mu\left(V_k \setminus \bigcup_{n=1}^{N} K_n\right) = \mu(X_N \setminus K_N) < \frac{1}{N}, \quad \forall N \in \mathbb{N},$$
showing that $\mu(V_k \setminus Y_k) = 0$.

Now denote
$$Y \doteq \bigcup_{k=1}^{\infty} Y_k \subset \mathcal{M}.$$

Then
$$\mu\left(Y^{\complement}\right) = \mu\left(\bigcap_{\ell=1}^{\infty} Y_\ell^{\complement}\right) = \mu\left(\bigcup_{k=1}^{\infty}\left[V_k \cap \bigcap_{\ell=1}^{\infty} Y_\ell^{\complement}\right]\right) \leq$$
$$\mu\left(\bigcup_{k=1}^{\infty} V_k \cap Y_k^{\complement}\right) = \mu\left(\bigcup_{\ell=1}^{\infty}(V_k \setminus Y_k)\right) \leq \sum_{k=1}^{\infty} \mu(V_k \setminus Y_k) = 0,$$
showing that $\mu$ is concentrated on $Y$. Since the cover $\{V_k\}_{k=1}^{\infty}$ is locally finite, the assertion of the proposition holds. $\square$

**Remark 2** *Proposition 1 above barely falls short of showing that a locally finite Borel measure on a Polish space is concentrated on a locally compact subspace; if countable disjoint unions of compact sets were always subsets of locally compact sets, their locally finite union would also be a subset of a locally compact set. But this is indeed not the case, which we will demonstrate on a simple example.*

**Example 1** *Let $\mathcal{M} = l^2(\mathbb{R})$ and let $\gamma : l_c^2(\mathbb{Q}) \to \mathbb{N}$ be any enumeration of the countable dense subset $l_c^2(\mathbb{Q}) \subset l^2(\mathbb{R})$ of rational sequences with finite support. For any Borel subset $U \subset l^2(\mathbb{R})$ set*
$$\mu(U) \doteq \sum_{x \in U \cap l_c^2(\mathbb{Q})} \frac{1}{\gamma(x)^2}.$$

*This defines a bounded Borel measure $\mu$ on the Polish space $l^2(\mathbb{R})$, which is indeed concentrated on the countable union $l_c^2(\mathbb{Q})$ of compact subsets (singletons). But there does not exist a locally compact subspace $Y$ such that $l_c^2(\mathbb{Q}) \subset Y \subset l^2(\mathbb{R})$. Indeed, if such a space existed then the point $0$ would have a compact neighbourhood $K \Subset Y$. Then also $K \Subset l^2(\mathbb{R})$ and thus $\overline{K}^{l^2} = K$. But $\exists \epsilon > 0$ such that $\mathcal{B}_\epsilon(0) \cap l_c^2(\mathbb{Q}) \subset K$, and hence,*
$$\overline{\mathcal{B}_\epsilon(0) \cap l_c^2(\mathbb{Q})}^{l^2} = \overline{\mathcal{B}_\epsilon(0)} \subset \overline{K}^{l^2} = K,$$
*which is a contradiction, because the closed balls are not compact in $l^2(\mathbb{R})$.*

Below is a slightly more sophisticated example of the same nature.

**Example 2** *Let $\mathcal{M} = l^2(\mathbb{R})$ and let*
$$\mathcal{M} = \bigsqcup_{k=1}^{\infty} \mathcal{M}_k, \quad \mathcal{M}_k \doteq \left\{x \in l^2(\mathbb{R}) \mid k-1 \leq \|x\|_2 < k\right\}, \quad \forall k \in \mathbb{N}.$$

*For every $k \in \mathbb{N}$ let $\gamma_k : l_c^2(\mathbb{Q}) \cap \mathcal{M}_k \to \mathbb{N}$ be any enumeration of points. Fix any non-negative sequence $\{a_k\}_{k=1}^{\infty} \in [0, +\infty)^{\mathbb{N}}$. For every Borel subset $U \subset l^2(\mathbb{R})$ set*
$$\mu(U) \doteq \sum_{x \in U \cap l^2(\mathbb{Q})} m(x), \quad m(x) \doteq \frac{a_k}{\gamma_k(x)^2}, \quad \forall x \in l^2(\mathbb{Q}) \cap \mathcal{M}_k, \quad \forall k \in \mathbb{N}.$$

Then $\mu$ is a Borel measure on $l^2(\mathbb{R})$, which is bounded or infinite depending on whether $\sum_k a_k$ converges or not. Weather $\mu$ is bounded or not, the measures of all balls satisfy $0 < \mu(\mathcal{B}_\epsilon(x)) < \infty$, $\forall x \in l^2(\mathbb{R})$, $\forall \epsilon > 0$. As in the example before, $\mu$ is concentrated on $l_c^2(\mathbb{Q})$ which is not a subset of any locally compact subspace.

On the other hand, the main assumptions of Theorem 1 - separability of $\mathcal{M}$ and inner regularity of $\mu$ - can be achieved at once if $\mathcal{M}$ is assumed $\sigma$-compact.

**Lemma 1** *If $\mathcal{M}$ is $\sigma$-compact, then it is separable, and every semifinite (see Definition 211F in [Fre01]) Borel measure on it is inner regular.*

**Proof:** Let $\mathcal{M} = \bigcup_{n=1}^{\infty} K_n$ with $K_n \Subset \mathcal{M}$ compact for $n \in \mathbb{N}$. Each $K_n$ is a compact metric space and hence Polish. Thus, the separability of $\mathcal{M}$ is obvious. Let

$$f \doteq \bigoplus_{n=1}^{\infty} \mathrm{id}_{K_n}, \quad f : \bigoplus_{n=1}^{\infty} K_n \to \mathcal{M},$$

then $f$ is continuous. The topological sum of countably many Polish spaces is Polish, and therefore $\mathcal{M}$ is a Souslin space. By Theorem 423E in [Fre03], all open subsets of $\mathcal{M}$ are $K$-analytic, and by Proposition 432C in [Fre03], $\mu$ is inner regular. $\square$

**Remark 3** *Suppose that the finite measure of balls (and hence local finiteness of $\mu$) is taken for granted. Then the $\sigma$-compactness of the space $\mathcal{M}$ guarantees the existence of a polar decomposition by Lemma 1 and Theorem 1. On the other hand, if $\mathcal{M}$ is separable and complete, then by Proposition 1 we can restrict $\mathcal{M}$ to a $\sigma$-compact subspace of full $\mu$-measure.*

## 3 The existence of polar decomposition

As noted before, a polar decomposition is a particular case of a measure disintegration. There are at least two sufficiently detailed modern expositions of measure disintegration that we know of: [Bog07.II] by Bogachev and [Fre03] by Fremlin. Theorem 10.4.8 together with Corollary 10.4.10 in [Bog07.II] have the merit of delivering a slightly stronger result, a regular conditional measure as opposed to a mere disintegration, which is valid also for signed measures. The deficiency of the treatment in [Bog07.II] is in the author's preoccupation with bounded measures only, which they merrily admit in the introduction to the first volume [Bog07.I]. Instead, Fremlin's section 452 in [Fre03] delivers a disintegration for $\sigma$-finite positive measures, which is just fine for our purposes. Otherwise the two expositions of the subject are essentially comparable.

The main result of this note is the following existence theorem.

**Theorem 1** *Let $(\mathcal{M}, d)$ be a separable metric space and $x_0 \in \mathcal{M}$. Let $\mu$ be a Borel measure on $(\mathcal{M}, \Sigma)$ such that:*

1. *$\mu$ is inner regular[1] (i.e., tight), that is,*

$$(\forall A \in \Sigma) \quad \mu(A) = \sup \{\mu(K) | \quad K \in \Sigma, \quad K \subset A, \quad K \text{ compact}\}.$$

---
[1] In fact, the assumption 1. can be relaxed to $\mu$ being countably compact in terms of Definition 451B in [Fre03].

2. Balls $\mathcal{B}_r(x_0) \subset \mathcal{M}$ of any radius[2] $r \in [0, R)$ have finite measure $\mu(\mathcal{B}_r(x_0)) < \infty$.

Then $(\mathcal{M}, d, \mu)$ admits a polar decomposition at $x_0$.

**Proof:** Recall the function $\rho_{x_0} \in C(\mathcal{M}, [0, R))$ defined as in (1) and the Borel measure $\nu_{x_0} = \mu \circ \rho_{x_0}^{-1}$ on $[0, R)$, so that $\rho_{x_0} : (\mathcal{M}, \mu) \to ([0, R), \nu_{x_0})$ is an inverse-measure-preserving map as in Definition 235G in [Fre03]. Since $[0, R)$ is Lindelöf, $\nu_{x_0}$ is $\sigma$-finite if and only if it is locally finite, and that is guaranteed by $\nu_{x_0}([0, r)) = \mu(\mathcal{B}_r(x_0)) < \infty$ for $\forall r > 0$. The countable generation of the Borel $\sigma$-algebra of $\mathcal{M}$ is given by the separability of $\mathcal{M}$. Now by Exercise 452X(l) in [Fre03], there exists a disintegration $\{\omega_r\}_{r \in [0, R)}$ of $\mu$ over $\nu_{x_0}$ consistent with $\rho_{x_0}$, which is essentially unique. For $\forall r \in [0, R)$, the object $\omega_r$ is a Borel probability measure on $\mathcal{M}$. Since $[0, R)$ is countably separated (Lemma 343E in [Fre02]), by Proposition 452G(c) in [Fre03] the disintegration $\{\omega_r\}_{r \in [0, R)}$ is strongly consistent with $\rho_{x_0}$ (see Definition 452E in [Fre03]). Thus, every $\omega_r$ is concentrated on $\rho_{x_0}^{-1}(\{r\}) = \mathcal{S}_r(x_0)$. By Proposition 452F in [Fre03], for every measurable $f : \mathcal{M} \to \mathbb{C}$ such that $\int_{\mathcal{M}} f(x) d\mu(x)$ makes sense, we have

$$\int_{\mathcal{M}} f(x) d\mu(x) = \int_{[0, R)} \int_{\mathcal{S}_r(x_0)} f(x) d\omega_r(x) d\nu_{x_0}(r).$$

The existence of a polar decomposition is established. $\square$

**Corollary 1** *Let $(\mathcal{M}, d)$ be a $\sigma$-compact metric space and $x_0 \in \mathcal{M}$. Let $\mu$ be a Borel measure on $\mathcal{M}$ such that all open balls centred at $x_0$ have finite $\mu$-measure. Then $(\mathcal{M}, d, \mu)$ admits a polar decomposition at $x_0$.*

**Proof:** Since all open balls $\mathcal{B}_r(x_0)$ have finite measure, $\mu$ is locally finite and hence semifinite. By Lemma 1, $\mathcal{M}$ is separable and $\mu$ is inner regular. It remains to apply Theorem 1. $\square$

For a metric space, separability, second countability and the Lindelöf property are mutually equivalent. Below we give a version of the above theorem which does not require the separability of $\mathcal{M}$, but instead puts harsher conditions on the measure $\mu$. Note that there is a discrepancy between the definition of a Radon measure space in [Fre03] and that of a Radon measure in most of the standard literature. Commonly, a Radon measure is a locally finite inner regular Borel measure. But Fremlin defines a Radon measure space what would normally be called a locally determined Borel measure space with a complete Radon measure (see Definition 411H(b) in [Fre03]).

**Theorem 2** *Let $(\mathcal{M}, d)$ be a metric space and $x_0 \in \mathcal{M}$. Let $\mu$ be a Borel measure on $(\mathcal{M}, \Sigma)$ such that:*

1. *$(\mathcal{M}, \Sigma, \mu)$ is a Radon measure space in the sense of Definition 411H(b) in [Fre03].*

2. *Balls $\mathcal{B}_r(x_0) \subset \mathcal{M}$ of any radius $r \in [0, R)$ have finite measure $\mu(\mathcal{B}_r(x_0)) < \infty$.*

*Then $(\mathcal{M}, d, \mu)$ admits a polar decomposition at $x_0$.*

**Proof:** Since the measure $\nu_{x_0}$ is still $\sigma$-finite, by Theorem 211L(c) in [Fre01] it is strictly localizable. By Proposition 452O in [Fre03], there exists a disintegration $\{\omega_r\}_{r \in [0, R)}$ of $\mu$ over $\nu_{x_0}$ consistent with $\rho_{x_0}$, and every $(\mathcal{M}, \omega_r)$ is a bounded Radon measure space. The rest proceeds as before. $\square$

---
[2]See the discussion around formula (3) regarding an economic choice of the bound $R$.

**Remark 4** *Note that Theorem 2 produces a polar decomposition where $(\mathcal{S}_r(x_0), \omega_r)$ is a bounded Radon measure space as per Definition 411H(b) in [Fre03].*

# 4 The absolute continuity of the measure $\nu_{x_0}$

To understand the nature of the question of absolute continuity of the measure $\nu_{x_0}$ in Definition 2 with respect to the Lebesgue measure, we will demonstrate very simple examples. The measure $\nu_{x_0}$ reflects the measure $\mu$ as distributed among spheres $\mathcal{S}_r(x_0)$, and its behaviour depends on the homogeneity of the measure $\mu$ as well as the commensurability of different spheres $\mathcal{S}_r(x_0)$. If $\mu$ is sufficiently homogeneous (whatever that means in a given context) then a non-absolutely continuous measure $\nu_{x_0}$ may arise due to an irregular foliation by non-commensurate spheres (again, whatever that means). On the plane, take $\mathcal{M}$ to be the union of an interval $[0, a] \times \{0\}$ on the horizontal axis and any arc on the unit circle centred at $(0, 0)$, and let $x_0 = (0, 0)$. Let $\mu$ be the Lebesgue measure on $\mathcal{M}$ measuring the length of curves in the usual way. Then the measure $\nu_{x_0}$ will have a singular part concentrated at $r = 1$. On the other hand, for any other $x_0 \in \mathcal{M}$ the measure $\nu_{x_0}$ would be absolutely continuous. But we can add arcs of various circles to produce more, infinitely many points $x_0 \in [0, a]$ with singular measures $\nu_{x_0}$.

Meanwhile, if the measure $\mu$ is not homogeneous enough then singular $\nu_{x_0}$ can arise even on ideally shaped metric spaces $\mathcal{M}$. Take $\mathcal{M}$ to be any standard smooth geometric object with Euclidean distance metric, e.g., the unit interval $[0, 1]$, and let $\mu$ contain a single point mass at any point. Then no matter where we take the point $x_0 \in [0, 1]$, the measure $\nu_{x_0}$ will not be absolutely continuous. The extreme example is the Cantor measure $\mu$, the Borel measure on $[0, 1]$ given by the Cantor function, which is singular everywhere on $[0, 1]$ [Fol99]. Here, too, the measure $\nu_{x_0}$ is purely singular for all $x_0 \in [0, 1]$.

Conditions on $(\mathcal{M}, d, \mu)$ and $x_0 \in \mathcal{M}$ under which $\nu_{x_0}$ is absolutely continuous are a delicate subject which should be studied separately. One thing that can be noticed is that the problems arising from the non-homogeneity of $\mu$ are far harder than those due to the uneven shape of $\mathcal{M}$.

# 5 Riemannian manifolds

Let us first consider the problem of finite measure of open balls in a Riemannian manifold.

**Remark 5** *Let $(\mathcal{M}, g)$ be a smooth connected Riemannian manifold with its geodesic distance $d_g$ and Riemannian volume $\upsilon_g$. For every ball $\mathcal{B}_r(x_0)$, if the Ricci curvature is bounded from below on $\mathcal{B}_r(x_0)$ then its Riemannian volume is finite, $\upsilon_g(\mathcal{B}_r(x_0)) < \infty$ .*

This statement follows immediately from the Gromov-Bishop-Günther comparison theorem, see Theorem 8.45 in [Gra04].

Note that the $\sigma$-compactness of a connected manifold is obvious, therefore Corollary 1 is applicable to Riemannian manifolds with Ricci curvature bounded from below on balls, and Borel measures $\mu$ having bounded Radon-Nikodym derivative with respect to the Riemannian volume. The latter condition guarantees finite $\mu$-measure of open balls, while their $\upsilon_g$-measures are already finite by the above remark.

Let us now turn to the absolute continuity of the measure $\nu_{x_0}$.

**Proposition 2** *Let $(\mathcal{M}, g)$ be a smooth connected Riemannian manifold with its geodesic distance $d_g$ and Riemannian volume $\upsilon_g$. Suppose that the Ricci curvature is bounded from below on every ball $\mathcal{B}_r(x_0)$. For every Borel measure $\mu$ on $\mathcal{M}$ that is absolutely continuous with respect to $\upsilon_g$, and for every point $x_0$, the measure $\nu_{x_0}$ is absolutely continuous.*

**Proof:** It suffices to prove the statement for $\mu = \upsilon_g$. Fix $x_0 \in \mathcal{M}$ and take $a \in (0, R)$, with $R$ as in (3). Since the Ricci curvature is bounded from below on $\mathcal{B}_a(x_0)$, by Gromov-Bishop-Günther comparison theorem (Theorem 8.45 in [Gra04]) there exists $c_a > 0$ such that

$$1 \geq \frac{\mu\left(\overline{\mathcal{B}_r(x_0)}\right)}{c_a r^n} \geq \frac{\mu\left(\overline{\mathcal{B}_s(x_0)}\right)}{c_a s^n}, \quad n \doteq \dim \mathcal{M},$$

for all $0 \leq r \leq s \leq a$. Thus,

$$0 \leq \mu\left(\overline{\mathcal{B}_s(x_0)}\right) - \mu\left(\overline{\mathcal{B}_r(x_0)}\right) \leq \frac{\mu\left(\overline{\mathcal{B}_r(x_0)}\right)}{r^n}(s^n - r^n) \leq nc_a a^{n-1}(s - r),$$

which shows that the function

$$r \mapsto \mu\left(\overline{\mathcal{B}_r(x_0)}\right) = \nu_{x_0}([0, r]) \tag{4}$$

is Lipschitz and hence absolutely continuous on any compact interval $[0, b] \subset [0, a)$. Because this is true for every $a \in (0, R)$, we find that the function (4) is absolutely continuous on every compact interval $[0, b] \subset [0, R)$. Therefore, by the fundamental theorem of Lebesgue integral calculus, the measure $\nu_{x_0}$ is absolutely continuous. □

A very similar argument can be found in Proposition 2.1 in [Skr19]. Thanks to Karen Hovhannisyan for pointing out this paper to us.

**Corollary 2** *Let $(\mathcal{M}, g)$ be a smooth connected Riemannian manifold with its geodesic distance $d_g$ and Riemannian volume $\upsilon_g$. Suppose that the Ricci curvature is bounded from below on every ball $\mathcal{B}_r(x_0)$. Let $\mu$ be a Borel measure on $\mathcal{M}$ such that $d\mu(x) = \dot{\mu}(x)d\upsilon_g(x)$ for $\dot{\mu} \in L^\infty(\mathcal{M}, \upsilon_g)$. Then at every $x_0 \in \mathcal{M}$, the metric measure space $(\mathcal{M}, d, \mu)$ admits a polar decomposition in the form given in Remark 1.*

**Remark 6** *The Ricci curvature is bounded globally on Riemannian manifolds with bounded geometry, including compact manifolds and homogeneous spaces. More generally, by Hopf-Rinow theorem, in a complete Riemannian manifold all balls are relatively compact, and thus the Ricci curvature is automatically bounded on each $\mathcal{B}_r(x_0)$. This makes the application of Corollary 2 completely straightforward in complete Riemannian manifolds.*

Note that for complete Riemannian manifolds taken with their natural volume measure, there are more geometric and explicit polar decompositions available. See, e.g., Proposition III.3.1 in [Chav06].

# Invariant sub-Riemannian structures on Lie groups

Let $G$ be a connected real Lie group, equipped with a left-invariant sub-Riemannian structure $(G, \mathcal{D}, \langle, \rangle)$ and the corresponding Carnot-Carathéodory metric $d$ [ABB19]. It is clear that $(G, d)$ is complete as a metric space (Corollary 7.51 in [ABB19]), so that by Proposition 3.47 in [ABB19], all balls $\mathcal{B}_r(x_0)$ are relatively compact, and hence have finite measure $\mu(\mathcal{B}_r(x_0)) < \infty$ for any locally finite Borel measure $\mu$. On the other hand, the $\sigma$-compactness of $G$ is also clear, and by Corollary 1 we establish the existence of a polar decomposition at every point $x_0 \in G$.

However, the question whether the measure $\nu_{x_0}$ is absolutely continuous appears to be open today.

# Declarations

## Funding

This work was supported by the FWO Odysseus 1 grant G.0H94.18N: Analysis and Partial Differential Equations and by the Methusalem programme of the Ghent University Special Research Fund (BOF) (Grant number 01M01021). M.R. was also supported by the EPSRC grants EP/R003025/2 and EP/V005529/1.

## Disclosure of potential conflicts of interest

Not applicable

## Compliance with ethical standards

Not applicable

## Data availability statement

Not applicable

## Code availability

Not applicable

# Ակնարկ մետրիկական տարածություններում բեռային վերլուծության մասին

Ժիրայր Ավետիսյան,* Մայքլ Ռուժանսկի†

September 1, 2023

**Սեղմագիր**

Էվկլիդյան տարածության բեռային կոորդինատների նմանակը՝ բեռային վերլուծությունը մետրիկական տարածությունում, երբ լավ սահմանված է, կարող է շատ օգտակար լինել բավարար ռեգուլյար չափի նկատմամբ ինտեգրալների հետ աշխատելիս։ Այս ակնարկում մենք դիտարկում ենք այդպիսի բեռային վերլուծությունների հետ կապված տեխնիկական մանրամասները։

**Հիմնաբառեր։** բեռային վերլուծություն, մետրիկական տարածություն, բեռային կոորդինատներ, Բորելի չափ, Ռադոնի չափ

## 1 Ներածություն

Թռչնի թռիչքի բարձրությունից դիտելով՝ չափելի ֆունկցիաների՝ անվերջ չափերի նկատմամբ ինտեգրալների հետ աշխատելիս հաճախ տեխնիկապես շատ հարմար է ամբողջական ինտեգրալը դիտարկել որպես կրկնակի ինտեգրալ, որտեղ չափի «անվերջությունն» ամբողջությամբ ընկնում է արտաքին ինտեգրալի վրա, ինչպես ստորև։

$$\int_{\mathcal{M}} f(x) d\mu(x) = \int_{[0,R)} \int_{\mathcal{S}_r} f(x) d\omega_r(x) d\nu(r), \quad R \in (0, +\infty]:$$

Այստեղ ենթադրվում է, որ կիրառված է չափային տարածությունների $(\mathcal{M}, \mu) \xrightarrow{\rho} ([0, R), \nu)$ տրցակ, $\mathcal{S}_r = \rho^{-1}(\{r\})$, ինչպես նաև

$$\mu(\cdot) = \int_{[0,R)} \omega_r(\cdot) d\nu(r)$$

տեսքի ապահնտեգրում (disintegration), որտեղ $\omega_r$-ը վերջավոր է և կենտրոնացված $\mathcal{S}_r$-ի վրա $\nu$-հ.ա. (համարյա ամեն՝ ըստ $\nu$-ի) $r \in [0, R)$-ի համար։ Եթե հնարավոր է $\mathcal{M}$-ը դիտարկել որպես $[0, R) \times \mathcal{S}$ արտադրյալ, այսինքն՝ գոյություն ունի Բորելի համարյա իզոմորֆիզմ $([0, R) \times \mathcal{S}, \tilde{\mu}) \xrightarrow{\Psi} \mathcal{M}$ այնպիսին, որ $\rho \circ \Psi = \pi_1$ (այստեղ $\pi_1(r, \varphi) = r$), ապա կարելի է ներմուծել $\mathcal{M} \ni x \mapsto (r, \varphi) \in [0, R) \times \mathcal{S}$ նոր փոփոխականներն

*Department of Mathematics: Analysis, Logic and Discrete Mathematics, Ghent University, Belgium, zhirayr.avetisyan@ugent.be

†Department of Mathematics: Analysis, Logic and Discrete Mathematics, Ghent University, Belgium, and School of Mathematical Sciences, Queen Mary University of London, United Kingdom. m.ruzhansky@ugent.be

այնպես, որ

$$\int_{\mathcal{M}} f(x)d\mu(x) = \int_{[0,R)} \int_{\mathcal{S}} \tilde{f}(r,\varphi)d\tilde{\omega}_r(\varphi)d\nu(r), \quad \tilde{f}(r,\varphi) = f \circ \Psi(r,\varphi), \quad d\tilde{\omega}_r(\varphi) = d\omega_r \circ \Psi(r,\varphi):$$

Այստեղ $\rho$ ֆունկցիայի ընտրությունը բավականին կամայական է, բայց մետրիկական տարածություններում կա դրա բնական թեկնածու՝ տրված $x_0 \in \mathcal{M}$ սկզբնակետից հեռավորությունը: Եթե $\mathcal{M} = \mathbb{R}^n$ and $d\mu(x) = \dot{\mu}(x)dx$, որտեղ $\dot{\mu} \in C^{\infty}(\mathbb{R}^n)$-ն ողորկ դրական խտություն է, ապա գործ ունենք պարզապես բևեռային կոորդինատների անցման հետ.

$$\int_{\mathbb{R}^n} f(x)d\mu(x) = \int_0^{\infty} \int_{\mathcal{S}_r} f(r,\varphi)\omega(r,\varphi)d\varphi dr, \quad \forall f \in L_{\text{loc}}^1(\mathbb{R}^n),$$

$$\mathcal{S}_r \doteq \{x \in \mathbb{R}^n \mid |x| = r\}, \quad \omega(r,\varphi) = \frac{d\mu(x(r,\varphi))}{d\varphi dr}:$$

Ընդհանուր դեպքում դիցուք $\mathcal{M} = (\mathcal{M}, d)$-ն մետրիկական տարածություն է, իսկ $\mu$-ն՝ չափ $\mathcal{M}$-ի $\Sigma$ Բորելի $\sigma$-հանրահաշվի վրա: Այս հոդվածում բոլոր չափերը կհամարվեն դրական, $\mu(\cdot) \geq 0$, և $\sigma$-ադդիտիվ՝

$$(\forall \{A_k\}_{k=1}^{\infty} \subset \Sigma) \quad (\forall k, l \in \mathbb{N}) \, k \neq l \Rightarrow A_k \cap A_l = \emptyset \quad \Rightarrow \quad \mu\left(\bigcup_{k=1}^{\infty} A_k\right) = \sum_{k=1}^{\infty} \mu(A_k):$$

Այսպիսի ցանկացած չափ մենք կանվանենք Բորելի չափ՝ չնայած, որ որոշ աղբյուրներ պահանջում են, որ Բորելի չափերն ունենան լրացուցիչ հարմար հատկություններ: Ընտրենք $x_0 \in \mathcal{M}$ սկզբնակետ և $\rho_{x_0} \in C(\mathcal{M}, [0, R))$-ով նշանակենք հետևյալ կերպ տրված ֆունկցիան.

$$\rho_{x_0}(x) = d(x, x_0), \quad \forall x \in \mathcal{M}: \tag{1}$$

Մենք այստեղ խոսելու ենք $x_0$ կենտրոնով և $r \in [0, R)$ շառավղով $\mathcal{S}_r(x_0)$ գնդոլորտների և $\mathcal{B}_r(x_0)$ գնդերի մասին.

$$\mathcal{S}_r(x_0) \doteq \rho_{x_0}^{-1}(\{r\}), \quad \mathcal{B}_r(x_0) \doteq \rho_{x_0}^{-1}([0,r)) = \bigcup_{0 \leq s < r} \mathcal{S}_s(x_0), \quad \forall r \in [0, R):$$

Բորելի $\nu_{x_0}$ չափը $[0, R)$-ի վրա սահմանված է այսպես՝

$$\nu_{x_0} \doteq \mu \circ \rho_{x_0}^{-1}, \tag{2}$$

և այդպիսով

$$\nu_{x_0}([0,r)) = \mu(\mathcal{B}_r(x_0)), \quad \forall r \in [0, R):$$

**Սահմանում 1** *Կասենք, որ մետրիկական չափային $(\mathcal{M}, d, \mu)$ տարածությունը թույլատրում է $x_0 \in \mathcal{M}$ կետում բևեռային վերլուծություն, եթե գոյություն ունի $\mathcal{M}$-ի վրա Բորելի վերջավոր $\omega_r$ չափերի $r \mapsto \omega_r$ դասը $\nu_{x_0}$-հ.ա. $r \in [0, R)$-ի համար այնպես, որ ցանկացած չափելի $f : \mathcal{M} \to \mathbb{C}$ ֆունկցիայի համար, որի*

$$\int_{\mathcal{M}} f(x)d\mu(x)$$

*ինտեգրալն իմաստալից է (վերջավոր կամ անվերջ), տեղի ունեն հետևյալ պայմանները.*

1. Հետնյալ արտապատկերումը չափելի է ըստ Բորելի՝

$$[0, R) \ni r \mapsto \int_{\mathcal{M}} f(x) d\omega_r(x)\text{:}$$

2. $\nu_{x_0}$-հ.ա. $r \in [0, R)$-ի համար $\omega_r$ չափը կենտրոնացած է $\mathcal{S}_r(x_0)$-ի վրա, այսինքն՝

$$\int_{\mathcal{M}} f(x) d\omega_r(x) = \int_{\mathcal{S}_r(x_0)} f(x) d\omega_r(x)\text{:}$$

3. Ունենք

$$\int_{\mathcal{M}} f(x) d\mu(x) = \int_{[0,R)} \int_{\mathcal{S}_r(x_0)} f(x) d\omega_r(x) d\nu_{x_0}(r) :$$

**Նկատառում 1** *Եթե $\nu_{x_0}$ չափը $[0, R)$-ի վրա բացարձակ անընդհատ է Լեբեգի չափի նկատմամբ, ապա գոյություն ունի ըստ Բորելի չափելի, ոչ բացասական, լոկալ ինտեգրելի $\dot{\nu}_{x_0}$ ֆունկցիա այնպես, որ $d\nu_{x_0}(r) = \dot{\nu}_{x_0}(r) dr$: Այդ դեպքում $\dot{\nu}_{x_0}$ արտադրիչը կարելի է ներառել $\omega_r$ չափերի մեջ այնպես, որ Սահմանում 1-ի մաս 3.-ը դառնա*

$$\int_{\mathcal{M}} f(x) d\mu(x) = \int_0^R \int_{\mathcal{S}_r(x_0)} f(x) d\omega_r(x) dr\text{:}$$

*Սա համապատասխանում է [RuVe19]-ում բևեռային վերլուծության (1.2) սահմանմանը. այս դեպքում այդ աշխատանքում կիրառված մասերով ինտեգրման մեթոդները կիրառելի են:*

Չափերի տեսության տեսանկյունից բևեռային վերլուծությունն ընդամենը չափի ապախնտեգրման մի դեպք է, և այս հոդվածում մեր նպատակն է ուսումնասիրել մետրիկական տարածության մեջ այս տեսակի չափի ապախնտեգրման համար անհրաժեշտ պայմանները: Այն հարցը, թե արդյոք $\nu_{x_0}$ չափը բացարձակ անընդհատ է Լեբեգի չափի նկատմամբ, առնչություն չունի ապախնտեգրման հետ, այլ արտացոլում է $\mathcal{M}$-ի երկրաչափական տեսքն ու $\mu$-ի բաշխումը:

## 2 Ենթադրությունների քննարկում

**Գնդերի վերջավոր չափը։** Սահմանում 1-ով տրված բևեռային վերլուծությունն օգտակար է, երբ պայման 3.-ը բավարարված է ոչ տրիվիալ կերպով: Եկեք արտագրենք այն չափերի լեզվով՝

$$\mu(\cdot) = \int_{[0,R)} \int_{\mathcal{S}_r(x_0)} \omega_r(\cdot) d\nu_{x_0}(r)\text{:}$$

Եթե $\mu$ չափը վերջավոր է, ապա ամեն ինչ պարզ է, բայց բևեռային վերլուծության ներմուծման իմաստը հաճախ հենց անվերջ չափի հետ աշխատանքն է: Եթե $\mu$ չափը լոկալ վերջավոր չէ, ապա անխուսափելիորեն $C(\mathcal{M}) \not\subset L^1_{\text{loc}}(\mathcal{M}, \mu)$, և անհրաժեշտ է դառնում աշխատել չափելի ֆունկցիաների շատ ավելի բազմությունների հետ: Կարելի է ասել, որ այդպիսի չափերը չունեն լայն կիրառություն ֆունկցիոնալ անալիզում, և ենթադրելով $\mu$ չափը լոկալ վերջավոր՝ ընդհանրության լուրջ խախտում տեղի չի ունենում:

Երբ $\mu$-ն լոկալ վերջավոր է $x_0$-ում, դիցուք

$$R_* \doteq \inf\{r \in [0, R) \mid \mu(\mathcal{B}_r(x_0)) = +\infty\} \in (0, +\infty], \tag{3}$$

պայմանով, որ $\inf \emptyset = R$։ Եթե $R_* < R$, ապա ամեն չափելի $A$ բազմության համար, որն ունի $\rho_{x_0}(A) \subset [R_*, R)$ պայմանով, ստանում ենք $\mu(A) = +\infty$, որն անօգուտ է դարձնում Սահմանում 1-ի 3. բանաձևը $r \geq R_*$ տիրույթում։ Ինչ վերաբերում է տրված $x_0$ կետի շուրջ բնեռային վերլուծությանը, կարելի է $\mathcal{M}$-ի փոխարեն դիտարկել $\mathcal{M}' \doteq \mathcal{B}_{R^*}(x_0)$ ենթատարածությունը և դնել $R \doteq R_*$։ Այս կերպ $r \in [0, R)$ շառավղով բոլոր $\mathcal{B}_r(x_0)$ գնդերը կունենան վերջավոր չափ։ Մասնավորապես $\mu$ չափը կլինի լոկալ վերջավոր $\mathcal{M}$-ի վրա։

**$(\mathcal{M}, d)$-ի սեպարաբելությունն ու լրիվությունը։** Մետրիկական $(\mathcal{M}, d)$ տարածության, և այդպիսով՝ նաև իր Բորելի $\sigma$-հանրահաշվի և իր վրա բոլոր Լեբեգի տարածությունների սեպարաբելության անշահ կարևոր է ֆունկցիոնալ անալիտիկ արդյունքների մեծ մասի համար։ Ուստի լուրջ սահմանափակում չի լինի ենթադրել, որ $(\mathcal{M}, d)$-ն սեպարաբել է։ Սակայն ստորև մենք կտեսնենք Թեորեմ 2-ում, որ բնեռային վերլուծություն կարելի է ստանալ նույնիսկ առանց սեպարաբելության, եթե պատրաստ ենք լրիվության և այլ պայմաններ դնել $\mu$-ի վրա։

Եթե մետրիկական $(\mathcal{M}, d)$ տարածությունը լրիվ չէ, ապա փակ գնդերը կարող են կոմպակտ չլինել և ունենալ անվերջ չափ, նույնիսկ երբ $\mu$-ն լոկալ վերջավոր է, իսկ $\mathcal{M}$-ը՝ լոկալ կոմպակտ։ Իրոք, այդպիսին է $\mathcal{M} = (0, 1)$-ն Էվկլիդյան մետրիկայով և $d\mu(x) = dx/x$ չափով։ $\mathcal{M}$-ի լրիվության ենթադրությունը կանխում է լոկալ վերջավոր Բորելի չափերի հնարավոր արտասովոր վարքը։ Թեռնս մեր հիմնական արդյունքները՝ Թեորեմ 1-ն ու Թեորեմ 2-ը, հիմնված չեն $(\mathcal{M}, d)$-ի լրիվության ենթադրության վրա։

Մյուս կողմից, եթե հաշվում ենք սեպարաբելության և լրիվության սահմանափակումների հետ, ապա հայտնվում ենք լեհական (Polish) մետրիկական $(\mathcal{M}, d)$ տարածությունների տիրույթում։ Լեհական մետրիկական տարածությունների վրա Բորելի չափերը շատ հնազանդ են և տեխնիկապես շատ հարմար։

**$\mathcal{M}$-ի $\sigma$-կոմպակտությունը։** Այս խորագրի ներքո կներկայացնենք երկու տեխնիկական պնդումներ, որոնք կարելի է համարել մեկը մյուսի մասնակի հակադարձը։ Միասին վերցրած՝ դրանք հուշում են, որ $\mathcal{M}$-ի $\sigma$-կոմպատկությունը շատ հեռու չէ բնեռային վերլուծությունների համար ճիշտ համատեքստը լինելուց՝ չնայած, որ այն խիստ անհրաժեշտ չէ։ Հիշենք, որ տոպոլոգիական տարածությունը կոչվում է $\sigma$-կոմպատկ, եթե այն հաշվելի հատ կոմպակտ ենթատարածությունների միավորում է։

Նախ ցույց տանք, որ եթե համաձայն ենք մետրիկական $(\mathcal{M}, d)$ տարածությունը համարել լեհական, ապա կարող ենք փոխարենը $(\mathcal{M}, d)$-ը համարել $\sigma$-կոմպակտ։

**Պնդում 1** *Լեհական տարածության վրա ամեն լոկալ վերջավոր Բորելի չափ կենտրոնացած է $\sigma$-կոմպակտ ենթատարածության վրա, որը որոշ հաշվելի, ճառավող կոմպակտ բազմությունների ընտանիքների լոկալ վերջավոր հաշվելի միավորումն է։*

**Ապացույց։** Դիցուք $\mathcal{M}$-ը լեհական տարածություն է, իսկ $\mu$-ն՝ լոկալ վերջավոր Բորելի չափ։ Դիցուք $\{U_x\}_{x \in \mathcal{M}}$-ը բաց բազմությունների այնպիսի ընտանիք է, որ $\forall x \in \mathcal{M}$-ի համար $\mu(U_x) < \infty$՝ ելնելով $\mu$-ի լոկալ վերջավոր հատկությունից։ Որպես մետրիզացվող տարածություն՝ $\mathcal{M}$-ը պարակոմպակտ է [Sto48], ուստի գոյություն ունի $\{U_x\}_{x \in \mathcal{M}}$-ի լոկալ վերջավոր բարելավում $\{V_\alpha\}_{\alpha \in A}$։ Բայց $\mathcal{M}$-ը նաև երկրորդ հաշվելի է (second countable) և այդպիսով՝ Լինդելոֆի, ուստի $\{V_\alpha\}_{\alpha \in A}$-ն ունի հաշվելի, լոկալ վերջավոր ենթածածկ՝ $\{V_k\}_{k=1}^\infty$։

Ընտրենք $k \in \mathbb{N}$: Պարզ է, որ $\exists x \in \mathcal{M}$ այնպիսին, որ $V_k \subset U_x$, և այդպիսով $\mu(V_k) \leq \mu(U_x) < \infty$: Որպես Սուսլինի $\mathcal{M}$ տարածության բաց ենթաբազմություն՝ $V_k$-ն ինքը Սուսլինի տարածություն է (Լեմմա 6.6.5, [Bog07.II]), և $\mu$-ի սահմանափակումը $V_k$-ի վրա Բորելի վերջավոր չափ է: Ըստ [Bog07.II]-ի Թեորեմ 7.4.3-ի՝ $\mu$-ն պիրկ (tight) է $V_k$-ի վրա, այսինքն՝

$$(\forall \epsilon > 0)(\exists K_\epsilon \Subset V_k \text{ կոմպակտ}) \, \mu(V_k \setminus K_\epsilon) < \epsilon:$$

Դիցուք $X_1 \doteq V_k$, իսկ կոմպակտ $K_1 \Subset X_1$ ընտրված է ըստ վերը գրվածի, և $\mu(X_1 \setminus K_1) < 1$: Ապա $X_2 \doteq X_1 \setminus K_1$-ը Սուսլինի տարածության բաց ենթաբազմություն է, ուստի՝ Սուսլինի տարածություն, և $\mu$-ն, սահմանափակված $X_2$-ի վրա, տալիս է Բորելի վերջավոր չափ, որը կրկին պիրկ է՝ ըստ նույն թեորեմի: Ինդուկտիվ շարունակելով այս գործընթացը՝ $n \in \mathbb{N}$-րդ քայլում ընտրում ենք կոմպակտ $K_n \Subset X_n$ այնպես, որ $\mu(X_n \setminus K_n) < \frac{1}{n}$: Նշանակենք

$$Y_k \doteq \bigcup_{n=1}^{\infty} K_n \subset X:$$

Ուրեմն

$$\mu(V_k \setminus Y_k) \leq \mu\left(V_k \setminus \bigcup_{n=1}^{N} K_n\right) = \mu(X_N \setminus K_N) < \frac{1}{N}, \quad \forall N \in \mathbb{N},$$

ինչը ցույց է տալիս, որ $\mu(V_k \setminus Y_k) = 0$:

Այժմ նշանակենք

$$Y \doteq \bigcup_{k=1}^{\infty} Y_k \subset \mathcal{M}:$$

Այսպիսով

$$\mu\left(Y^{\complement}\right) = \mu\left(\bigcap_{\ell=1}^{\infty} Y_\ell^{\complement}\right) = \mu\left(\bigcup_{k=1}^{\infty} \left[V_k \cap \bigcap_{\ell=1}^{\infty} Y_\ell^{\complement}\right]\right) \leq$$

$$\mu\left(\bigcup_{k=1}^{\infty} V_k \cap Y_k^{\complement}\right) = \mu\left(\bigcup_{\ell=1}^{\infty}(V_k \setminus Y_k)\right) \leq \sum_{k=1}^{\infty} \mu(V_k \setminus Y_k) = 0,$$

ինչը ցույց է տալիս, որ $\mu$-ն կենտրոնացած է $Y$-ի վրա. Քանի որ $\{V_k\}_{k=1}^{\infty}$ ծածկը լոկալ վերջավոր է, պնդումն ապացուցված է: $\square$

**Նկատառում 2** *Վերոնշյալ Պնդում 1-ը, կարելի էր կարծել, թե ցույց է տալիս, որ լոկալ վերջավոր Բորելի չափը լեհական տարածության վրա կենտրոնացած է լոկալ կոմպակտ ենթատարածության վրա. եթե հաշվելի հատ չհատվող կոմպակտ բազմությունների միավորումը միշտ լիներ լոկալ կոմպակտ բազմության ենթաբազմություն, այդպիսի միավորումների լոկալ վերջավոր միավորումը նույնպես կլիներ լոկալ կոմպակտ բազմության ենթաբազմություն: Բայց դա իրոք այպես չէ, և մենք դա կցուցադրենք մի պարզ օրինակով:*

**Օրինակ 1** *Դիցուք $\mathcal{M} = l^2(\mathbb{R})$, և $\gamma : l_c^2(\mathbb{Q}) \to \mathbb{N}$-ով նշանակենք վերջավոր կրիչով (support) ռացիոնալ հաջորդականությունների՝ հաշվելի և խիտ $l_c^2(\mathbb{Q}) \subset l^2(\mathbb{R})$ ենթաբազմության որևէ համարակալում: Ցանկացած $U \subset l^2(\mathbb{R})$ Բորելի ենթաբազմության համար նշանակենք*

$$\mu(U) \doteq \sum_{x \in U \cap l_c^2(\mathbb{Q})} \frac{1}{\gamma(x)^2}:$$

*Այսպես սահմանվում է Բորելի վերջավոր $\mu$ չափը լեհական $l^2(\mathbb{R})$ տարածության վրա, որն իրոք կենտրոնացած է կոմպակտ (միակետ) բազմությունների հաշվելի միավորման՝ $l_c^2(\mathbb{Q})$-ի վրա: Բայց գոյություն չունի $Y$ լոկալ*

կոմպակտ ենթատարածություն այնպիսին, որ $l_c^2(\mathbb{Q}) \subset Y \subset l^2(\mathbb{R})$։ Հիրավի, եթե այդպիսի տարածություն գոյություն ունենար, ապա 0 կետը կունենար մի կոմպակտ շրջակայք՝ $K \in Y$։ Այդ դեպքում նաև $K \in l^2(\mathbb{R})$, ուստի նաև $\overline{K}^{l^2} = K$։ Բայց $\exists \epsilon > 0$ այնպիսին, որ $\mathcal{B}_\epsilon(0) \cap l_c^2(\mathbb{Q}) \subset K$, և այդպիսով

$$\overline{\mathcal{B}_\epsilon(0) \cap l_c^2(\mathbb{Q})}^{l^2} = \overline{\mathcal{B}_\epsilon(0)} \subset \overline{K}^{l^2} = K,$$

ինչը հակասություն է, քանի որ $l^2(\mathbb{R})$-ում փակ գնդերը կոմպակտ չեն։

Ստորև բերենք նույն բնույթի մի փոքր ավելի բարդ օրինակ։

**Օրինակ 2** Դիցուք $\mathcal{M} = l^2(\mathbb{R})$, և նշանակենք

$$\mathcal{M} = \bigsqcup_{k=1}^\infty \mathcal{M}_k, \quad \mathcal{M}_k \doteq \left\{ x \in l^2(\mathbb{R}) \,\big|\, k-1 \leq \|x\|_2 < k \right\}, \quad \forall k \in \mathbb{N}։$$

Ամեն $k \in \mathbb{N}$-ի համար $\gamma_k : l_c^2(\mathbb{Q}) \cap \mathcal{M}_k \to \mathbb{N}$-ով նշանակենք կետերի որևէ համարակալում։ Ընտրենք որևէ ոչ բացասական $\{a_k\}_{k=1}^\infty \in [0, +\infty)^\mathbb{N}$ հաջորդականություն։ Ամեն Բորելի $U \subset l^2(\mathbb{R})$ ենթաբազմության համար նշանակենք

$$\mu(U) \doteq \sum_{x \in U \cap l^2(\mathbb{Q})} m(x), \quad m(x) \doteq \frac{a_k}{\gamma_k(x)^2}, \quad \forall x \in l^2(\mathbb{Q}) \cap \mathcal{M}_k, \quad \forall k \in \mathbb{N}։$$

Այսպիսով $\mu$-ն Բորելի չափ է $l^2(\mathbb{R})$-ի վրա, որը վերջավոր է կամ անվերջ՝ կախված նրանից, թե արդյոք $\sum_k a_k$-ն զուգամիտում է, թե ոչ։ Լինի $\mu$-ն վերջավոր, թե անվերջ, բոլոր գնդերի չափերը բավարարում են $0 < \mu(\mathcal{B}_\epsilon(x)) < \infty$ պայմանին, $\forall x \in l^2(\mathbb{R}), \forall \epsilon > 0$։ Ինչպես և առաջին օրինակում, $\mu$-ն կենտրոնացած է $l_c^2(\mathbb{Q})$-ի վրա, որը որևէ լոկալ կոմպակտ ենթատարածության ենթաբազմություն չէ։

Մյուս կողմից՝ Թեորեմ 1-ի հիմնական ենթադրությունները՝ $\mathcal{M}$-ի սեպարաբելությունն ու $\mu$-ի պիրկությունը, ստացվում են միանգամից, եթե $\mathcal{M}$-ը ենթադրում ենք $\sigma$-կոմպակտ։

**Լեմմա 1** Եթե $\mathcal{M}$-ը $\sigma$-կոմպակտ է, ապա այն սեպարաբել է, և ամեն կիսավերջավոր (տես Սահմանում 211F, [Fre01]) Բորելի չափ նրա վրա պիրկ է (նույնիսկը՝ ներքին ռեգուլյար)։

**Ապացույց:** Դիցուք $\mathcal{M} = \bigcup_{n=1}^\infty K_n$, որտեղ $K_n \Subset \mathcal{M}$-ը կոմպակտ է բոլոր $n \in \mathbb{N}$-երի համար։ Ամեն $K_n$-ը կոմպակտ մետրիկական, ուստի և լեհական տարածություն է։ Այդպիսով $\mathcal{M}$-ի սեպարաբելությունն ակնհայտ է։ Եթե նշանակենք

$$f \doteq \bigoplus_{n=1}^\infty \mathrm{id}_{K_n}, \quad f : \bigoplus_{n=1}^\infty K_n \to \mathcal{M},$$

ապա $f$-ն անընդհատ է։ Հաշվելի հատ լեհական տարածությունների ուղիղ գումարը լեհական է, ուստի $\mathcal{M}$-ը Սուսլինի տարածություն է։ Ըստ [Fre03]-ի Թեորեմ 423E-ի՝ $\mathcal{M}$-ում բոլոր բաց բազմությունները $K$-անալիտիկ են, ուստի, ըստ [Fre03]-ի Պնդում 432C-ի, $\mu$-ն պիրկ է։ $\square$

**Նկատառում 3** Համարենք, որ բոլոր գնդերի վերջավոր չափը (և այդպիսով $\mu$-ի լոկալ վերջավոր լինելը) ապահովված է։ Այդ դեպքում $\mathcal{M}$ տարածության $\sigma$-կոմպակտությունը երաշխավորում է բնեռային վերլուծության գոյությունը Լեմմա 1-ի և Թեորեմ 1-ի օգնությամբ։ Մյուս կողմից, եթե $\mathcal{M}$-ը սեպարաբել է ու լրիվ, ապա Պնդում 1-ի օգնությամբ մենք կարող ենք $\mathcal{M}$-ի փոխարեն վերցնել լրիվ $\mu$-չափ ունեցող $\sigma$-կոմպակտ ենթատարածություն։

# 3 Բնեռային վերլուծության գոյությունը

Ինչպես նշեց վերը, բնեռային վերլուծությունը չափի ապախնտեգրման մասնավոր դեպք է: Մեզ հայտնի են չափի ապախնտեգրման առնվազն երկու բավարար մանրամասն նկարագրություններ. [Bog07.II]-ը՝ Բոգաչովի հեղինակությամբ, և [Fre03]-ը՝ Ֆրեմլինի: [Bog07.II]-ի Թեորեմ 10.4.8-ն ու Հետևանք 10.4.10-ն ունեն այն առավելությունը, որ փոքր-ինչ ավելի ուժեղ արդյունք են տալիս՝ ռեգուլյար պայմանական չափ, ոչ թե պարզապես ապախնտեգրում, և գործում են նաև նշանավոր չափերի համար: Այս աշխատանքի թերությունը հեղինակի հակումն է դեպի զուտ վերջավոր չափեր, որը նա խաղաղությամբ խոստովանում է առաջին հատորի ներածականում [Bog07.I]: Փոխարենը՝ Ֆրեմլինի բաժին 452-ը [Fre03]-ում տրամադրում է դրական $\sigma$-վերջավոր չափերի ապախնտեգրում, ինչը բավարար է մեր նպատակների համար: Այլապես թեմային այս երկու անդրադարձները համեմատելի են:

Այս ակնարկի գլխավոր արդյունքը հետևյալ գոյության թեորեմն է:

**Թեորեմ 1** *Դիցուք $(\mathcal{M}, d)$-ն սեպարաբել մետրիկական տարածություն է, և $x_0 \in \mathcal{M}$: Դիցուք $\mu$-ն այնպիսի Բորելի չափ է $(\mathcal{M}, \Sigma)$-ի վրա, որ.*

1. *$\mu$-ն ներքին ռեգուլյար է՝ (այսինքն՝ պիրկ), ասել է, թե*

$$(\forall A \in \Sigma)(\forall \epsilon > 0)(\exists K_\epsilon \subset A) \quad K_\epsilon \text{ կոմպակտ} \quad \wedge \quad \mu(A \setminus K_\epsilon) < \epsilon:$$

2. *Կամայական $r \in [0, R)$ շառավղով[2] գնդերն ունեն վերջավոր չափ՝ $\mu(\mathcal{B}_r(x_0)) < \infty$:*

*Այս պարագայում գոյություն ունի $\mu$-ի բնեռային վերլուծություն $\mathcal{M}$-ում $x_0$-ի շուրջ:*

**Ապացույց:** Հիշենք (1)-ում սահմանված $\rho_{x_0} \in C(\mathcal{M}, [0, R))$ ֆունկցիան ու $\nu_{x_0} = \mu \circ \rho_{x_0}^{-1}$ Բորելի չափը $[0, R)$-ի վրա, և նկատենք, որ $\rho_{x_0} : (\mathcal{M}, \mu) \to ([0, R), \nu_{x_0})$-ն հակադարձ չափապահպան (inverse-measure-preserving) արտապատկերում է՝ ըստ [Fre03]-ի Սահմանում 235G-ի: Քանի որ $[0, R)$-ը Լինդելոֆի է, $\nu_{x_0}$-ն $\sigma$-վերջավոր է այն և միայն այն դեպքում, երբ այն լոկալ վերջավոր է, իսկ դա երաշխավորված է $\nu_{x_0}([0, r)) = \mu(\mathcal{B}_r(x_0)) < \infty$ պայմանով $\forall r > 0$-ի համար: $\mathcal{M}$-ի Բորելի $\sigma$-հանրահաշվի հաշվելի ծնումը հետևում է $\mathcal{M}$-ի սեպարաբելությունից: Այժմ [Fre03]-ի Վարժություն 452X(l)-ից հետևում է, որ գոյություն ունի $\mu$-ի ապախնտեգրում $\nu_{x_0}$-ի նկատմամբ՝ $\{\omega_r\}_{r \in [0,R)}$, որը համատեղելի է $\rho_{x_0}$-ի հետ և էապես միակն է: $\forall r \in [0, R)$-ի համար $\omega_r$ օբյեկտը Բորելի հավանականային չափ է $\mathcal{M}$-ի վրա: Քանի որ $[0, R)$-ը հաշվելիորեն տարանջատված է (countably separated, տես Լեմմա 343E,[Fre02]), ըստ [Fre03]-ի Պնդում 452G(c)-ի՝ $\{\omega_r\}_{r \in [0,R)}$ ապախնտեգրումը խիստ համատեղելի է $\rho_{x_0}$-ի հետ (տես Սահմանում 452E,[Fre03]): Այսպիսով յուրաքանչյուր $\omega_r$ կենտրոնացած է $\rho_{x_0}^{-1}(\{r\}) = \mathcal{S}_r(x_0)$-ի վրա: Ըստ [Fre03]-ի Պնդում 452F-ի՝ ամեն չափելի $f : \mathcal{M} \to \mathbb{C}$-ի համար, որի $\int_\mathcal{M} f(x) d\mu(x)$-ը իմաստալից է, ունենք

$$\int_\mathcal{M} f(x) d\mu(x) = \int_{[0,R)} \int_{\mathcal{S}_r(x_0)} f(x) d\omega_r(x) d\nu_{x_0}(r):$$

Այսպիսով բնեռային վերլուծության գոյությունը հաստատված է: $\square$

---
[1]Իրականում ենթադրություն 1.-ը կարող է թուլացվել այնպես, որ $\mu$-ն լինի հաշվելի կոմպակտ [Fre03]-ի Սահմանում 451B-ի եզրաբանությամբ, կամ ունենա կոմպակտ մոտարկող դաս ինչպես [Bog07.II]-ի Սահմանում 1.4.6-ում:

[2]Տես բանաձև (3)-ի շուրջ քննարկումը $R$ սահմանի խսայողական ընտրության մասին:

**Հետևանք 1** *Դիցուք $(\mathcal{M}, d)$-ն $\sigma$-կոմպակտ մետրիկական տարածություն է, և $x_0 \in \mathcal{M}$: Դիցուք $\mu$-ն այնպիսի Բորելի չափ է $\mathcal{M}$-ի վրա, որ $x_0$ կենտրոնով բոլոր բաց գնդերն ունեն վերջավոր $\mu$-չափ: Այդ դեպքում գոյություն ունի $\mu$-ի բևեռային վերլուծություն $\mathcal{M}$-ում $x_0$-ի շուրջ:*

**Ապացույց:** Քանի որ բոլոր $\mathcal{B}_r(x_0)$ բաց գնդերը վերջավոր չափ ունեն, $\mu$-ն լոկալ վերջավոր է, և այդպիսով նաև՝ կիսավերջավոր: Ըստ Լեմմա 1-ի՝ $\mathcal{M}$-ը սեպարաբել է, իսկ $\mu$-ն՝ պիրկ: Մնում է միայն կիրառել Թեորեմ 1-ը: □

Մետրիկական տարածության համար սեպարաբելությունը, երկրորդ հաշվելիությունն ու Լինդելոֆի հատկությունը փոխադարձ համարժեք են: Ստորև ներկայացնում ենք վերը տրված թեորեմի մի տարբերակ, որը չի պահանջում $\mathcal{M}$-ի սեպարաբելությունը, բայց փոխարենն ավելի խիստ պայմաններ է դնում $\mu$ չափի վրա: Նկատենք, որ կա անհամապատասխանություն մի կողմից [Fre03]-ում Ռադոնի չափային տարածության սահմանման, իսկ մյուս կողմից ստանդարտ գրականության մեջ Ռադոնի չափի սահմանման միջև: Սովորաբար Ռադոնի չափ համարվում է լոկալ վերջավոր և ներքին ռեգուլյար (պիրկ) Բորելի չափը: Բայց Ֆրեմլինը Ռադոնի չափային տարածություն է անվանում այն, ինչ մնացյալ գրականության մեջ կկոչվեր լոկալ որոշյալ Բորելի չափային տարածություն լրիվ Ռադոնի չափով (տես Սահմանում 411H(b),[Fre03]):

**Թեորեմ 2** *Դիցուք $(\mathcal{M}, d)$-ն մետրիկական տարածություն է, և $x_0 \in \mathcal{M}$: Դիցուք $\mu$-ն այնպիսի Բորելի չափ է $(\mathcal{M}, \Sigma)$-ի վրա, որ․*

1. *$(\mathcal{M}, \Sigma, \mu)$-ն Ռադոնի չափային տարածություն է [Fre03]-ի Սահմանում 411H(b)-ի իմաստով:*

2. *Կամայական $r \in [0, R)$ շառավղով գնդերն ունեն վերջավոր չափ՝ $\mu(\mathcal{B}_r(x_0)) < \infty$:*

*Այս պարագայում գոյություն ունի $\mu$-ի բևեռային վերլուծություն $\mathcal{M}$-ում $x_0$-ի շուրջ:*

**Ապացույց:** Քանի որ $\nu_{x_0}$ չափը կրկին $\sigma$-վերջավոր է, ըստ [Fre01]-ի Թեորեմ 211L(c)-ի՝ այն խիստ տեղայնացվող է (strictly localizable): Ըստ [Fre03]-ի Պնդում 452O-ի՝ գոյություն ունի $\mu$-ի ապահնտեգրում $\nu_{x_0}$-ի նկատմամբ՝ $\{\omega_r\}_{r \in [0,R)}$, որը համատեղելի է $\rho_{x_0}$-ի հետ, և յուրաքանչյուր $(\mathcal{M}, \omega_r)$ Ռադոնի վերջավոր չափային տարածություն է: Մնացյալը՝ նախորդ թեորեմի համանմանությամբ: □

**Նկատառում 4** *Նկատենք, որ Թեորեմ 2-ը տրամադրում է բևեռային վերլուծություն, որտեղ $(\mathcal{S}_r(x_0), \omega_r)$-ը Ռադոնի վերջավոր չափային տարածություն է՝ ըստ [Fre03]-ի Սահմանում 411H(b)-ի:*

# 4 $\nu_{x_0}$ չափի բացարձակ անընդհատությունը

Սահմանում 2-ում $\nu_{x_0}$ չափի՝ Լեբեգի չափի նկատմամբ բացարձակ անընդհատության խնդրի բնույթը հասկանալու համար այստեղ դիտարկենք շատ պարզ օրինակներ: $\nu_{x_0}$ չափն արտացոլում է $\mu$ չափը՝ բաժշկված ըստ $\mathcal{S}_r(x_0)$ գնդոլորտների, և նրա վարքը կախված է $\mu$ չափի համասեռությունից, ինչպես նաև տարբեր $\mathcal{S}_r(x_0)$ գնդոլորտների համեմատելիությունից: Եթե $\mu$-ն բավարար համասեռ է (որևից իմաստով), ապա ոչ բացարձակ անընդհատ $\nu_{x_0}$ կարող է առաջանալ անհամեմատելի գնդոլորտներով անկանոն շերտավորման պատճառով (կրկին, որևից իմաստով): Հարթության վրա դիցուք $\mathcal{M}$-ը հորիզոնական առանցքի վրա $[0, a] \times \{0\}$ հատվածի և $(0, 0)$ կենտրոնով միավոր շրջանի որեն ադերի միավորումն է, իսկ $x_0 = (0, 0)$: Ենթադրենք $\mu$-ն Լեբեգի չափն է $\mathcal{M}$-ի վրա, որը չափում է կորի երկարությունը սովորական կերպով: Այստեղ $\nu_{x_0}$ չափը ունի եզակի (singular) մաս՝ կենտրոնացած $r = 1$ կետում: Սյուս կողմից՝ ցանկացած այլ $x_0 \in \mathcal{M}$ կետի համար $\nu_{x_0}$ կլիներ բացարձակ անընդհատ։ Ավելացնելով տարբեր շոշանների ադեղներ՝ կարող ենք առաջացնել ավելի շատ, անվերջ քանակի $x_0 \in [0, a]$ կետեր, որոնց $\nu_{x_0}$ չափերը եզակի են։

Միաժամանակ, եթե $\mu$ չափը բավարար համասեռ չէ, ապա եզակի $\nu_{x_0}$ կարող է առաջանալ նույնիկ կատարյալ տեսքով մետրիկական $\mathcal{M}$ տարածությունում։ Կարելի է վերցնել որպես $\mathcal{M}$ ցանկացած սովորական ողորկ երկրաչափական մարմին Էվկլիդյան հեռավորությամբ, օրինակ՝ $[0, 1]$ միավոր հատվածը, իսկ $\mu$-ին թույլ տալ ունենալ ընդամենը մեկ կետային զանգված։ Այդ դեպքում ինչպես էլ ընտրենք $x_0 \in [0, 1]$ կետը, $\nu_{x_0}$ չափը չի լինի բացարձակ անընդհատ։ Ցայտահեղ օրինակ է Կանտորի $\mu$ չափը՝ Կանտորի ֆունկցիայով $[0, 1]$-ի վրա տրվող Բորելի չափը, որը եզակի է ամենուրեք $[0, 1]$-ի վրա [Fol99]։ Այս դեպքում նույնպես $\nu_{x_0}$-ն մաքուր եզակի է բոլոր $x_0 \in [0, 1]$-ների համար։

Պայմանները $(\mathcal{M}, d, \mu)$-ի և $x_0 \in \mathcal{M}$-ի վրա, որոնց դեպքում $\nu_{x_0}$-ն բացարձակ անընդհատ է, նուրբ խնդիր են ներկայացնում, որը պետք է առանձին ուսումնասիրել։ Թերևս կարելի է նկատել, որ $\mu$-ի անհամասեռությունից բխող խնդիրները շատ ավելի բարդ են, քան $\mathcal{M}$-ի ձևից կախված խնդիրները։

## 5   Ռիմանյան բազմաձևություններ

Եկեք նախ դիտարկենք Ռիմանյան բազմաձևության մեջ բաց գնդերի վերջավոր չափի խնդիրը։

**Նկատառում 5** *Դիցուք $(\mathcal{M}, g)$-ն ողորկ, կապակցված Ռիմանյան բազմաձևություն է իր գեոդեզիական $d_g$ հեռավորությամբ և Ռիմանյան $v_g$ ծավալով։ Յուրաքանչյուր $\mathcal{B}_r(x_0)$ գնդի համար, եթե Ռիչիի կորությունը ներքևից սահմանափակ է $\mathcal{B}_r(x_0)$-ի վրա, ապա գնդի Ռիմանյան ծավալը վերջավոր է՝ $v_g(\mathcal{B}_r(x_0)) < \infty$։*

Այս պնդումն անմիջապես հետևում է Գրոմով-Բիշոփ-Գյունթերի համեմատության թեորեմից, տես [Gra04]-ի Թեորեմ 8.45-ը։

Նկատենք, որ կապակցված բազմաձևության $\sigma$-կոմպակտությունն ակնհայտ է, ուստի Հետևանք 1-ը կիրառելի է այնպիսի Ռիմանյան բազմաձևությունների համար, որոնց Ռիչիի կորությունը ներքևից սահմանափակ է գնդերի վրա, և Բորելի այնպիսի $\mu$ չափերի համար, որոնք ունեն սահմանափակ Ռադոն-Նիկոդիմի ածանցյալ Ռիմանյան ծավալի նկատմամբ։ Այս վերջին պայմանը երաշխավորում է բաց գնդերի վերջավոր $\mu$-չափը, իսկ դրանց $v_g$-չափն արդեն իսկ վերջավոր է ըստ վերոնշյալ նկատառման։

Այժմ անդրադառնանք $\nu_{x_0}$ չափի բացարձակ անընդհատությանը։

**Պնդում 2** *Դիցուք $(\mathcal{M}, g)$-ն ողորկ, կապակցված Ռիմանյան բազմաձևություն է իր գեոդեզիական $d_g$ հեռավորությամբ և Ռիմանյան $v_g$ ծավալով։ Ենթադրենք, որ Ռիչիի կորությունը ներքևից սահմանափակ է ամեն $\mathcal{B}_r(x_0)$ գնդի վրա։ $\mathcal{M}$-ի վրա Բորելի յուրաքանչյուր $\mu$ չափի համար, որը բացարձակ անընդհատ է $v_g$-ի նկատմամբ, և ամեն $x_0$ կետի համար, $\nu_{x_0}$ չափը բացարձակ անընդհատ է։*

**Ապացույց:** Բավարար է պնդումն ապացուցել $\mu = v_g$-ի համար։ Ընտրենք $x_0 \in \mathcal{M}$ և $a \in (0, R)$, որտեղ $R$-ն՝ ինչպես (3)-ում։ Քանի որ Ռիչիի կորությունը ներքևից սահմանափակ է $\mathcal{B}_a(x_0)$-ի վրա, ըստ Գրոմով-Բիշոփ-Գյունթերի թեորեմի ([Gra04]-ի Թեորեմ 8.45) գոյություն ունի $c_a > 0$ այնպիսին, որ

$$1 \geq \frac{\mu\left(\overline{\mathcal{B}_r(x_0)}\right)}{c_a r^n} \geq \frac{\mu\left(\overline{\mathcal{B}_s(x_0)}\right)}{c_a s^n}, \quad n \doteq \dim \mathcal{M},$$

բոլոր $0 \leq r \leq s \leq a$-ի համար։ Այսպիսով՝

$$0 \leq \mu\left(\overline{\mathcal{B}_s(x_0)}\right) - \mu\left(\overline{\mathcal{B}_r(x_0)}\right) \leq \frac{\mu\left(\overline{\mathcal{B}_r(x_0)}\right)}{r^n}(s^n - r^n) \leq nc_a a^{n-1}(s-r),$$

ինչը ցույց է տալիս, որ

$$r \mapsto \mu\left(\overline{\mathcal{B}_r(x_0)}\right) = \nu_{x_0}([0,r]) \tag{4}$$

ֆունկցիան Լիպշիցյան է, ուստի՝ բացարձակ անընդհատ ցանկացած կոմպակտ $[0,b] \subset [0,a)$ միջակայքում։ Քանի որ սա ճշմարիտ է ցանկացած $a \in (0, R)$-ի համար, ստանում ենք, որ (4) ֆունկցիան բացարձակ անընդհատ է յուրաքանչյուր կոմպակտ $[0,b] \subset [0,R)$ միջակայքում։ Հետևաբար, ըստ Լեբեգի ինտեգրալ հաշվի հիմնարար թեորեմի, $\nu_{x_0}$ չափը բացարձակ անընդհատ է։ □

Բավականին համանման փաստարկ կիրառվել է [Skr19]-ի Պնդում 2.1-ում։ Շնորհակալություն Կարեն Հովհաննիսյանին այդ հոդվածը մեր ուշադրությանը ներկայացնելու համար։

**Հետևանք 2** *Դիցուք $(\mathcal{M}, g)$-ն ողորկ, կապակցված Ռիմանյան բազմաձևություն է իր գեոդեզիական $d_g$ հեռավորությամբ և Ռիմանյան $v_g$ ծավալով։ Ենթադրենք, որ Ռիչիի կորությունը ներքևից սահմանափակ է ամեն $\mathcal{B}_r(x_0)$ գնդի վրա։ Դիցուք $\mu$-ն այնպիսի Բորելի չափ է $\mathcal{M}$-ի վրա, որ $d\mu(x) = \dot{\mu}(x)dv_g(x)$ որևէ $\dot{\mu} \in L^\infty(\mathcal{M}, v_g)$-ի համար։ Այդ դեպքում յուրաքանչյուր $x_0 \in \mathcal{M}$ կետում մետրիկական չափային $(\mathcal{M}, d, \mu)$ տարածությունը թույլ է տալիս բնեռային վերլուծություն Նկատառում 1-ում տրված տեսքով։*

**Նկատառում 6** *Ռիչիի կորությունը գլոբալ սահմանափակ է սահմանափակ երկրաչափությամբ Ռիմանյան բազմաձևությունների վրա, ինչպիսիք են կոմպակտ բազմաձևություններն ու համասեռ տարածությունները։ Ընդհանրապես, ըստ Հոպֆ-Ռինովի թեորեմի, Ռիմանյան լրիվ բազմաձևության բոլոր գնդերը հարաբերական կոմպակտ են, և այդպիսով Ռիչիի կորությունն ինքնըստինքյան սահմանափակ է ամեն $\mathcal{B}_r(x_0)$-ի վրա։ Այս փաստի շնորհիվ Հետևանք 2-ի կիրառումը Ռիմանյան լրիվ բազմաձևություններում միանգամայն պարզ է։*

Նշենք, որ իրենց բնական ծավալային չափի հետ դիտարկվող Ռիմանյան լրիվ բազմաձևությունների համար գոյություն ունեն ավելի երկրաչափական և բացահայտ բնեռային վերլուծություններ։ Տես, օրինակ, [Cha06]-ի Պնդում III.3.1-ը։

# Ինվարիանտ Ենթաիմանյան կառուցվածքներ Լիի խմբերի վրա

Դիցուք $G$-ն կապակցված իրական Լիի խումբ է՝ կախավորված ձախ ինվարիանտ Ենթաիմանյան $(G, \mathcal{D}, \langle, \rangle)$ կառուցվածքով և համապատասխան Կառնո-Կարաթեոդորիի $d$ մետրիկայով [ABB19]։ Պարզ է, որ $(G, d)$-ն որպես մետրիկական տարածություն լրիվ է ([ABB19]-ի Հետևանք 7.51), ուստի, ըստ [ABB19]-ի Պնդում 3.47-ի, բոլոր $\mathcal{B}_r(x_0)$ գնդերը հարաբերական կոմպակտ են և այդպիսով ունեն վերջավոր $\mu(\mathcal{B}_r(x_0)) < \infty$ չափ Բորելի ցանկացած լոկալ վերջավոր $\mu$ չափի պարագայում։ Մյուս կողմից՝ $G$-ի $\sigma$-կոմպատկության նույնպես պարզ է, և Հետևանք 1-ի օգնությամբ կարող ենք հաստատել բնեռային վերլուծության գոյությունը ցանկացած $x_0 \in G$ կետում։

Սակայն այն հարցը, թե արդյոք $\nu_{x_0}$ չափը բացարձակ անընդհատ է, ըստ ամենայնի այսոր մնում է բաց։

### Շահերի հնարավոր բախման բացահայտում

Կիրառելի չէ

### Համապատասխանություն էթիկայի նորմերին

Կիրառելի չէ

### Տվյալների հասանելիության հայտարարություն

Կիրառելի չէ

### Կոդի հասանելիությունը

Կիրառելի չէ

## Գրականություն